\documentclass[11pt,twoside,reqno]{amsart}
\usepackage{import}
\usepackage{newclude}
\usepackage{anysize}
\usepackage{amsmath}
\usepackage{amsthm}
\usepackage{amsmath,amscd}
\usepackage[utf8]{inputenc}
\usepackage{amssymb}
\usepackage{stmaryrd}
\usepackage{wasysym}
\usepackage{mathrsfs}
\usepackage{hyperref}
\usepackage{cleveref}
\usepackage{enumitem}
\usepackage{dsfont}
\usepackage{soul}
\hypersetup{colorlinks=true,linkcolor=black,citecolor=black}
\usepackage{color}
\usepackage[all]{xy}
\usepackage{float}
\usepackage{bm}
\usepackage{mathtools}
\usepackage{thmtools}
\usepackage{setspace}
\usepackage{comment}
\usepackage{tikz}
\usepackage{tikz-cd}
\usepackage{mathdots}
\usepackage{graphicx}
\usepackage[myheadings]{fullpage}

\numberwithin{equation}{section}
\newcommand\mtop{.95in}
\newcommand\mbottom{.95in}
\newcommand\mleft{1in}
\newcommand\mright{1in}
\usepackage[top = \mtop, bottom = \mbottom, left = \mleft, right=\mright]{geometry}
\DeclareMathOperator{\val}{val}
\DeclareMathOperator{\Mat}{Mat}

\newtheorem{thm}{Theorem}[section]
\newtheorem{example}[thm]{Example}
\newtheorem{prop}[thm]{Proposition}

\newtheorem{lemma}[thm]{Lemma}

\newtheorem{cor}[thm]{Corollary}
\theoremstyle{definition}
\newtheorem{defi}[thm]{Definition}

\newtheorem{rmk}[thm]{Remark}

\newcommand\reallywidehat[1]{%
\savestack{\tmpbox}{\stretchto{%
  \scaleto{%
    \scalerel*[\widthof{\ensuremath{#1}}]{\kern-.6pt\bigwedge\kern-.6pt}%
    {\rule[-\textheight/2]{1ex}{\textheight}}
  }{\textheight}%
}{0.5ex}}%
\stackon[1pt]{#1}{\tmpbox}%
}
\DeclareSymbolFont{bbold}{U}{bbold}{m}{n}
\DeclareSymbolFontAlphabet{\mathbbold}{bbold}

\makeatletter
\def\@tocline#1#2#3#4#5#6#7{\relax
  \ifnum #1>\c@tocdepth 
  \else
    \par \addpenalty\@secpenalty\addvspace{#2}%
    \begingroup \hyphenpenalty\@M
    \@ifempty{#4}{%
      \@tempdima\csname r@tocindent\number#1\endcsname\relax
    }{%
      \@tempdima#4\relax
    }%
    \parindent\z@ \leftskip#3\relax \advance\leftskip\@tempdima\relax
    \rightskip\@pnumwidth plus4em \parfillskip-\@pnumwidth
    #5\leavevmode\hskip-\@tempdima
      \ifcase #1
       \or\or \hskip 1em \or \hskip 2em \else \hskip 3em \fi%
      #6\nobreak\relax
    \hfill\hbox to\@pnumwidth{\@tocpagenum{#7}}\par
    \nobreak
    \endgroup
  \fi}
\makeatother

\newcommand{\R}{\mathbb{R}}
\newcommand{\Z}{\mathbb{Z}}
\newcommand{\Q}{\mathbb{Q}}

\newcommand{\C}{\mathbb{C}}
\newcommand{\F}{\mathbb{F}}

\newcommand{\E}{\mathbb{E}}

\renewcommand{\l}{\lambda}

\renewcommand{\L}{\Lambda}

\newcommand{\G}{\mathbb{G}}

\newcommand{\bx}{\mathbf{x}}

\DeclareMathOperator{\GL}{GL}

\DeclareMathOperator{\SO}{SO}

\DeclareMathOperator{\GSp}{GSp}

\DeclareMathOperator{\Sig}{Sig}

\DeclareMathOperator{\GT}{GT}

\DeclareMathOperator{\Cov}{Cov}

\DeclareMathOperator{\SN}{SN}
\DeclareMathOperator{\diag}{diag}

\title{Gaussian universality of $p$-adic random matrix products via corners}
\author{Jiahe Shen}

\date{\today}
\begin{document}

\thanks{I thank my advisor Ivan Corwin for helpful advice on revisions and for providing me funding support with his NSF grant DMS-2246576 and Simons Investigator grant 929852; and Roger Van Peski, for his hints that inspire my research, suggestions on the materials, and valuable encouragement and comments for polishing my rhetorical skills.  
}

\maketitle

\begin{abstract} 
We establish the universality of the singular numbers in random matrix products over $\GL_n(\Q_p)$ as the number of products approaches infinity, with a fixed $n\ge 1$. We demonstrate that, under a broad class of distributions, which we term as ``split", the asymptotics of matrix products align with the sum of the singular numbers of the matrix corners. Specifically, when matrices are independent and identically distributed, we derive the strong law of large numbers and the central limit theorem. Our approach is inspired by Van Peski’s work \cite{Roger}, which examines products of $n\times n$ corners of Haar-distributed elements $\GL_N(\Z_p)$. We extend the method so that the criterion now works as long as the measures are left- and right-invariant under the multiplication of $\GL_n(\Z_p)$. Building on this approach for $\GL_n$, we further demonstrate similar universality for $\GSp_{2n}$ and discuss potential extensions to general split reductive groups. 
\end{abstract}

\textbf{Keywords: }\keywords{$p$-adic random matrices, universality, strong law of large numbers, central limit theorem, split reductive group}

\textbf{Mathematics Subject Classification (2020): }\subjclass{15B52 (primary); 15B33, 60B20 (secondary)}

\tableofcontents

\section{Introduction}

\subsection{Results.}
Random matrix theory, which explores the properties of matrix-valued random variables, is a vibrant and expanding area of research. A common approach in this field is to diagonalize the matrices and focus on their diagonal entries. For instance, given $n\le m$ and nonsingular complex matrix $A\in\Mat_{n\times m}(\C)$, there exists $U\in U(n),V\in U(m)$ such that $UAV=\diag_{n\times m}(\alpha_1,\ldots,\alpha_n)$, where $0<\alpha_1\le\ldots\le\alpha_n$. Motivated by the asymptotic results of the sum of random variables, such as the law of large numbers and central limit theorem, Bellman \cite{bellman1954limit}, followed by Furstenberg-Kesten \cite{furstenberg1960products}, studied the limits and fluctuations of the singular values in the product sequence $A_1,A_1A_2,\ldots$, where $A_1,A_2\ldots$ are random matrices in the real or complex field.

As a $p$-adic analogy of the above, given $n\le m$ and nonsingular $p$-adic matrix $A\in\Mat_{n\times m}(\Q_p)$, there exist matrices $U\in\GL_n(\Z_p),V\in\GL_m(\Z_p)$ such that $A=U\diag_{n\times m}(p^{\lambda_1},\ldots,p^{\lambda_n})V$, where $(\lambda_1\ge\ldots\ge\lambda_n)\in\Z^n$. The integers $\l_i$ are referred as the \textbf{singular numbers} of $A$, denoted $\SN(A) = (\l_1,\ldots,\l_n) = \l$. Furthermore, the \textbf{weight} of $\lambda$ is defined as $|\lambda|=\lambda_1+\ldots+\lambda_n$. In an analogy with the archimedean case, it is natural to study the asymptotic of the sequence $\SN(A_1),\SN(A_1A_2),\SN(A_1A_2A_3),\ldots$, where $A_1,A_2,A_3\ldots$ are independent random matrices from $\GL_n(\mathbb{Q}_p)$. Van Peski’s work \cite{Roger}, which inspired our method, establishes the strong law of large numbers and the central limit theorem for this sequence. In his framework, for every $k\ge 1$, the matrix $A_k$ is $n\times n$ corner of a Haar distributed matrix. Our work generalizes this result to a wide range (see ``split" in \Cref{defi: split}) of measures as long as they are invariant under the left-and right-multiplication of $\GL_n(\Z_p)$. In particular, when the random matrices $A_1,A_2,\ldots$ are i.i.d. (the $n^2$ entries of $A_1$ need not be i.i.d.), we prove the following theorem:

\begin{thm}\label{cor: i.i.d.}
Let $A_1,A_2,\ldots\in\GL_n(\Q_p)$ be i.i.d. random matrices whose distributions are invariant under both left- and right-multiplication by $\GL_n(\Z_p)$. Let $\lambda(k)=(\lambda_1(k),\ldots,\lambda_n(k))=\SN(A_1\ldots A_k)$. For every $1\le i\le n$, let $A_1^{(i)}\in\Mat_{(n-i+1)\times n}(\Q_p)$ denote the submatrix of $A_1$ consisting of the last $n-i+1$ rows and all $n$ columns. 
\begin{enumerate}[left=0pt]
\item (Strong law of large numbers) Suppose that for all $1\le i\le n$, the expectation $\E|\SN(A_1^{(i)})|<\infty$ exists. Then we have
\begin{align}(\frac{\l_1(k)}{k},\ldots,\frac{\l_n(k)}{k})\rightarrow &(\E|\SN(A_1^{(1)})|-\E|\SN(A_1^{(2)})|,\ldots,\nonumber\\
&\E|\SN(A_1^{(n-1)})|-\E|\SN(A_1^{(n)})|,\E|\SN(A_1^{(n)})|)\nonumber\end{align}
almost surely as $k$ goes to infinity. 

\item(Central limit theorem) In addition to the conditions above, assume that for all $1\le i\le n$, the expectation $\E|\SN(A_1^{(i)})^2|<\infty$ exists. Then we have weak convergence
$$\frac{(\l_1(k)-k\E|\SN(A_1^{(1)})|+k\E|\SN(A_1^{(2)})|,\ldots,\l_n(k)-k\E|\SN(A_1^{(n)})|)}{\sqrt{k}}\Rightarrow N(0,L\Sigma L^T)$$
as $k$ goes to infinity, where $L^T\in\Mat_{n\times n}(\Z_p)$ is the transpose of the matrix
$$L=\begin{pmatrix}1 & -1 & & &\\ & 1 & -1 & &\\ & & \ddots & \ddots &\\ & & & 1 & -1\\ & & & & 1\end{pmatrix}$$
and $\Sigma=\Cov_{1\le i,j\le n}(|\SN(A_1^{(i)})|,|\SN(A_1^{(j)})|)$ is the covariance matrix.
\end{enumerate}
\end{thm}

It is worth noting that the i.i.d. case only serves as a specific example of our method: it also applies to non-i.i.d. setting such as the result from Van Peski \cite{Roger} (see \eqref{item: Roger setting} of \Cref{example: split}). A different way of generalization is to replace $\GL_n$ by other split reductive groups, for instance the $\GSp_{2n}$ case in \Cref{thm: i.i.d for GSp}. As far as we know, this is the first demonstration of the equivalence between random matrix products and the sum of random matrix corners, a connection that we believe encompasses the vast majority of relevant research scenarios.

Before we delve into the rigorous proof, it’s helpful to outline the heuristic that motivates our method. Suppose $A_1,\ldots,A_k$ are already given, and we write $\lambda(k)=(\lambda_1(k),\ldots,\lambda_n(k))=\SN(A_1\cdots A_k)$. We are interested in determining the distribution of $\SN(A_1\cdots A_{k+1})$. While this matrix product may seem challenging to handle directly, we take advantage of the fact that the distribution is invariant under left- and right-multiplication by $\GL_n(\Z_p)$. This allows us to diagonalize the matrix $A_1\cdots A_k$ without changing the distribution and then analyze the singular numbers of the resulting matrix.

Thus, we examine the singular numbers
$$\lambda(k+1)=(\lambda_1(k+1),\ldots,\lambda_n(k+1))=\SN(\diag_{n\times n}(p^{\lambda_1(k)},\ldots,p^{\lambda_n(k)})A_{k+1}).$$
As $k$ goes to infinity, we expect to see the distance between the parts of $\lambda(k)$ also goes to infinity, i.e., $\lambda_1(k)$ is much larger than $\lambda_2(k)$, $\lambda_2(k)$ is much larger than $\lambda_3(k)$, and so on. Recall the well-known property that singular numbers can be derived from the minimal value of the determinant of minors (see \Cref{thm:submatrices_suffice}). Hence, it is reasonable to speculate that for all $1\le i\le n$, when we search for the $(n-i+1)\times(n-i+1)$ submatrix of $\diag_{n\times n}(p^{\lambda_1(k)},\ldots,p^{\lambda_n(k)})A_{k+1}$ whose determinant has the smallest value, the last $n-i+1$ rows must be picked. Therefore, our heuristic surmises that except for finitely many $k$, we have

$$\lambda_i(k+1)+\ldots+\lambda_n(k+1)-(\lambda_i(k)+\ldots+\lambda_n(k))=|\SN(A_k^{(i)})|,\quad\forall 1\le i\le n$$
where $A_k^{(i)}\in\Mat_{(n-i+1)\times n}(\Q_p)$ is the submatrix of $A_k$ consisting of the last $n-i+1$ rows and all $n$ columns. Hence, one would expect the following to share the same asymptotic:

\begin{enumerate}[left=0pt]
\item \label{item: matrix product}The sequence $\lambda(1),\lambda(2),\ldots$, where the $n$-dimension joint variable $\lambda(k)=(\lambda_1(k),\ldots,\lambda_n(k))\in\Sig_n:=\{\lambda=(\lambda_1,\ldots,\lambda_n)\in\Z^n\mid\lambda_1\ge\ldots\ge\lambda_n\}$ is the singular numbers of matrix product $A_1\cdots A_k$;
\item \label{item: add corner} The sequence $v(1),v(2),\ldots$, where the $n$-dimension joint variable $v(k)=(v_1(k),\ldots,v_n(k))\in\Z^n$ satisfies (different from $\lambda(k)\in\Sig_n$ which must be a signature, the non-increasing condition $v_1(k)\ge\ldots\ge v_n(k)$ might not hold, so $v(k)$ might not be an element of $\Sig_n$) 

$$v_i(k)+\ldots+v_n(k)=\sum_{j=1}^k |\SN(A_j^{(i)})|.$$
\end{enumerate} 

\eqref{item: add corner} is much easier to deal with than \eqref{item: matrix product}, since it simply adds independent random variables and no matrix product is involved. Nevertheless, this idea is not rigorous, and extra conditions are required (which will be named ``split" in \Cref{Proof}) to ensure the distance between the singular numbers $\lambda_1(k),\ldots,\lambda_n(k)$ would increase to infinity when $k$ goes to infinity. Our first step is to fill these gaps and give the precise statement in \Cref{thm: main theorem}; For the second step, we turn to the i.i.d case of \Cref{cor: i.i.d.} and prove the inequality \eqref{eq: inequality of corners} by linear algebra. Hence the sequence $v$ is split almost surely, and the strong law of large numbers and central limit theorem for i.i.d. matrix products are directly inherited from the sum of singular numbers of i.i.d corners.

\begin{rmk}
In Van Peski's original paper \cite{Roger}, there were also two sequences: the ``interacting dynamic'' $\lambda(0)=(0[n]),\lambda(1),\ldots$ that comes from the random matrix product, and the ``non-interacting dynamic'' $v(0)=(0[n]),v(1),\ldots$, which is the sum of independent random variables. Similar to our result,  \cite[Proposition 5.4]{Roger} shows that the difference between the sequences $\lambda$ and $v$ is almost surely bounded and therefore simplifies the problem into the sum of independent random variables. However, the two sequences $\lambda$ and $v$ are generated in a different way: first, Van Peski regards multiplying the matrix $A_k$ as the ascending dynamics in the explicit Hall-Littlewood process, and thus gives the joint distribution of products $\SN(A_1),\SN(A_1A_2),\ldots$ in the form of Hall-Littlewood polynomials (see the appendix for backgrounds.). This is feasible because the Truncated Haar ensemble in \cite[Theorem 1.3]{Roger} ensures the singular numbers $\SN(A_k)$ to be distributed with respect to the Hall-Littlewood measure; then, the two sequences $\lambda$ and $v$ are raised by relating them to an interacting particle system similar to PushTASEP. Our work, however, is purely linear algebra, which only concerns the corners coming from the above heuristic. This skips the Hall-Littlewood ascending process step and directly uses the matrix corners to generate the sequences $\lambda$ and $v$, which proceeds beyond the restriction of Hall-Littlewood distribution. 
\end{rmk}

\begin{rmk}\label{rmk: reductive group}
In the heuristic mentioned above, the only advantage of having distributions that are invariant under left- and right-multiplication by $\GL_n(\Z_p)$ is that it enables us to diagonalize the matrix $A_1\cdots A_k$ before multiplying by $A_{k+1}$. Therefore, we expect that the heuristic remains valid for a wide range of split reductive groups. In \Cref{reductive group}, we will state and prove the strong law of large numbers and central limit theorem for the $\GSp_{2n}$ case (\Cref{thm: i.i.d for GSp}), which turns out to have a similar description with the $\GL_n$ case in \Cref{cor: i.i.d.}. More generally, for any split reductive group $\G$, let $G=\G(\Q_p)$, and $K=\G(\Z_p)$ be the maximal compact subgroup of $G$. By Cartan decomposition in \cite[Section 4]{bruhat1972groupes}, there exist elements $U,V\in K$ such that $A=Ud(\lambda)V$, where $d(\lambda)$ is a normal form related to the cocharacter $\lambda\in P^+$ in the root system. We still denote $\SN(A)=\lambda$ as the singular numbers of $A$ and study the asymptotic behavior of the sequence $\SN(A_1),\SN(A_1A_2),\ldots,$ where $A_1,A_2,\ldots$ are independently chosen from $G$ with distributions invariant under both left- and right-multiplication by $K$. Further discussion on this part is postponed to \Cref{reductive group}.
\end{rmk}

Let us turn our attention to previous progress in the archimedean case. The study of asymptotic behavior in random matrices has attracted considerable attention from various fields, including statistical physics and theoretical computer science. Many works have been published on the asymptotics of random matrix products, including the contributions of Avila-Eskin-Viana \cite{avila2023continuitylyapunovexponentsrandom}, Newman \cite{newman1986distribution}, Vanneste \cite{vanneste2010estimating}, Forrester \cite{forrester2015asymptotics}, Lima-Rahibe \cite{lima1994exact}, Akemann-Ipsen \cite{akemann2015recent}, Liu-Wang-Wang \cite{liu2018lyapunov}. These works often focus on the understanding of \textbf{Lyapunov exponents} from both mathematical and physical perspectives. For random complex matrices $A_1,A_2,\ldots$, the $i^{th}$ Lyapunov exponent is defined as
$$\lim_{k\rightarrow\infty}\frac{1}{k}\log(\text{$i^{th}$ largest singular number of $A_1A_2\cdots A_k$})$$
Following the suggestions of \cite{Roger}, for the non-archimedean case, it is reasonable to regard the $i^{th}$ Lyapunov exponent for random matrices $A_1,A_2,\ldots\in\GL_n(\Q_p)$ as the limit
$$\lim_{k\rightarrow\infty}\frac{\lambda_{n-i+1}(k)}{k}=\lim_{k\rightarrow\infty}\frac{\SN_{n-i+1}(A_1\cdots A_k)}{k}$$
Thus, the strong law of large numbers in \Cref{cor: i.i.d.} gives the $i^{th}$ Lyapunov exponent for i.i.d matrices over $\Q_p$. 

On the other hand, various techniques have been developed to analyze the $p$-adic case, for which the settings may have little difference. Except from \cite{Roger}, Van Peski's work \cite{vanpeski2024padicdysonbrownianmotion} also concerns the i.i.d case as mentioned in \Cref{cor: i.i.d.}, but mainly focused on the case when $\SN(A_1)=(1,0,\ldots,0)$ is fixed. In his settings, the frequency of multiplying new random matrices follows the Poisson process. Also, if we view the product as random walks on groups with an ergodic theory perspective, one can find universality by generalizing Oseledets' multiplicative ergodic theorem. This is what Brofferio-Schapira \cite{MR2837130} did in their work, which proves a law of large numbers for products of i.i.d. random matrices over $\GL_n(\Q_p)$. The measures in \cite{MR2837130} are not invariant under the left- and right-multiplication of $\GL_n(\Z_p)$, but our result covers a much wider range of distributions for $\SN(A_k),k=1,2,\ldots$, and provides a more refined estimate such as the central limit theorem, which goes beyond the method in \cite{MR2837130}. Chhaibi \cite{chhaibi2017non} studied the probabilistic approach based on Jacquet's Whittaker function over a reductive group over $\GL_n(\Q_p)$. 

\begin{rmk}\label{rmk:any_local_field}
Throughout the paper, the matrices are over the $p$-adic field $\Q_p$ because it is the most commonly considered. Nevertheless, the results and proofs in this paper still hold for matrices over any non-archimedean local field $F$ with finite residue field, i.e., any finite algebraic extensions of $\Q_p$ or $\F_p((t))$. Let $\mathfrak{o}$ be the ring of integers of $F$, $\mathfrak{p}$ the maximal ideal of $\mathfrak{o}$, $\pi$ a generator of $\mathfrak{p}$, and $q=|\mathfrak{o}/\mathfrak{p}|$ be the order of the residue field. Since the Hecke algebra structure holds in this generality (see Chapter 5 of Macdonald's book), the only thing we need to do is to replace $\Q_p$ by $F$, $\Z_p$ by $\mathfrak{o}$, $p$ by $\pi$ (in the context of matrix entries), and setting $t=1/q$ in the Hall-Littlewood specializations in \Cref{HL and RMT}.
\end{rmk}

\subsection{Plan of paper.} In \Cref{Preliminaries}, we discuss $p$-adic random matrices and probabilistic constructions over them; in \Cref{Proof}, we give the strict proof of the heuristic above (\Cref{thm: main theorem}) and use this to prove the i.i.d. case in \Cref{cor: i.i.d.}; in \Cref{reductive group}, we apply our heuristic to random matrix product over arbitrary split reductive group. In particular, we prove the strong law of large numbers and central limit theorem for the general symplectic group $\GSp_{2n}$ case; In \Cref{HL and RMT}, we introduce the Hall-Littlewood polynomials and their connections to $p$-adic random matrices, thereby providing an alternating proof of \Cref{cor: i.i.d.}.

\section{Preliminaries}\label{Preliminaries}

In this section, we briefly summarize properties of random matrices over the $p$-adic field. 

\begin{defi}
Let $\Sig_n:=\{\lambda=(\lambda_1,\ldots,\lambda_n)\in\Z^n\mid\lambda_1\ge\ldots\ge\lambda_n\}$ denote the set of integer signatures of length $n$, which are non-increasing $n$-tuples of integers, and let $\Sig_n^+\subset\Sig_n$ denote set of signatures with all parts nonnegative. Given $\l = (\l_1,\ldots,\l_n) \in \Sig_n$, we refer to the integers $\l_i$ as the \textbf{parts} of $\l$. We set $|\l| := \sum_{i=1}^n \l_i,n(\lambda):=\sum_{i=1}^n (i-1)\lambda_i$, and $m_k(\l) = \#\{i\mid \l_i = k\}$. For $\l \in \Sig_n$ and $\mu \in \Sig_{n-1}$, write $\mu \prec_P \l$ if $\l_i \geq \mu_i$ and $\mu_i \geq \l_{i+1}$ for $1 \leq i \leq n-1$. We write $c[k]$ for the signature $(c,\ldots,c)$ of length $k$.
\end{defi}

Let us introduce the background of $p$-adic fields, which is a typical example of non-archimedean fields. Fix a prime $p$. Any nonzero rational number $r \in \Q^\times$ can be written as 
$$r=p^k\frac{a}{b},\quad a,b,k\in\Z, (a,p)=(b,p)=1.$$ 
Define $|\cdot|: \Q \to \R$ by setting $|r|_p = p^{-k}$ for $r$ as above, and $|0|_p=0$. Then $|\cdot|_p$ gives a norm on $\Q$ and $d_p(x,y) :=|x-y|_p$ gives a metric. We also denote $\val_p(r)=k$ for $r=p^k\frac{a}{b}$ as above and $\val_p(0) = \infty$, so $|r|_p = p^{-\val_p(r)}$. 

We define the \textbf{$p$-adic field} $\Q_p$ to be the set
$$\Q_p:=\{r=\sum_{i=k}^\infty a_ip^i\Big| k\in\Z,a_i\in\{0,1,\ldots,p-1\}\}$$
which is completion of $\Q$ with respect to the matrix $|\cdot|$. Also, the \textbf{$p$-adic integers} $\Z_p$ is defined to be the below subring of $\Q_p$:
$$\Z_p:=\{r\in\Q_p\Big|\; |r|\le 1\}=\{r=\sum_{i=0}^\infty a_ip^i\Big| a_i\in\{0,1,\ldots,p-1\}\}.$$

$\Q_p$ is equipped with an additive invariant Haar measure, which is unique if we require $\Z_p$ to have measure $1$. The restriction of this measure to $\Z_p$ is the unique Haar probability measure on $\Z_p$, whose pushforward under any quotient $r_n:\Z_p \to \Z_p/p^n\Z_p$ is the uniform probability measure. Therefore, the Haar probability measure over $\Z_p$ has the explicit interpretation that each $a_i$ is i.i.d. uniformly random from $\{0,\ldots,p-1\}$.

Similarly, $\GL_N(\Q_p)$ has a unique left- and right-multiplicative invariant measure for which $\GL_N(\Z_p)$ has measure $1$. The restriction of this measure to $\GL_N(\Z_p)$, which we denote by $M_{Haar}(\GL_N(\Z_p))$, pushes forward to $\GL_N(\Z/p^n\Z)$ and is the uniform measures on the finite groups $\GL_N(\Z/p^n\Z)$.

The following proposition comes from Cartan decomposition, which gives the orbits of nonsingular matrices in $\Mat_{n \times m}(\Q_p)$ under the multiplication of $\GL_n(\Z_p)\times\GL_m(\Z_p)$ on both sides.

\begin{prop}\label{prop:smith}
Let $n \leq m$. For any nonsingular $A \in \Mat_{n \times m}(\Q_p)$, there exist $U \in \GL_n(\Z_p), V \in \GL_m(\Z_p)$ such that $A = U\diag_{n \times m}(p^{\lambda_1},\ldots,p^{\lambda_n})V$ where $\lambda=(\lambda_1,\ldots,\lambda_n)\in \Sig_n$ is unique. 
\end{prop}

The parts $\lambda_i$ are known as the \textbf{singular numbers} of $A$, and we will also write $\SN(A) = \lambda = (\lambda_1,\ldots,\lambda_n)$. We also use the $\diag$ symbol for block matrices. From now on, unless otherwise specified, all the matrices are assumed to be nonsingular. The restriction of the Haar measure on $\GL_N(\Q_p)$ to the double coset 
$$\GL_N(\Z_p) \diag_{N\times N}(p^{\lambda_1},\ldots,p^{\lambda_N}) \GL_N(\Z_p),$$
normalized to be a distribution, is the unique $\GL_N(\Z_p) \times \GL_N(\Z_p)$-invariant distribution on the subset of $\GL_N(\Q_p)$ with singular numbers $\lambda$. These distributions are equivalently written as $U \diag_{N\times N}(p^{\lambda_1},\ldots,p^{\lambda_N}) V$ where $U,V$ are independently distributed by the Haar probability measure on $\GL_N(\Z_p)$. More generally, if $n \leq m$ and $U \in \GL_n(\Z_p), V \in \GL_m(\Z_p)$ are Haar distributed and $\mu \in \Sig_n$, then $U \diag_{n \times m}(p^\mu) V$ is invariant under $\GL_n(\Z_p) \times \GL_m(\Z_p)$ acting on the left and right, and is the unique such bi-invariant measure with singular numbers given by $\mu$. As a corollary, we have

\begin{cor}\label{cor: invariant measure}
Suppose $n\le m$, and $A\in\Mat_{n\times m}(\Q_p)$ is random, nonsingular with distribution invariant under the action of $\GL_n(\Z_p)\times\GL_m(\Z_p)$. In other words, $A$ is an element in the probability space $(\Omega,\mathcal{F},G)$, where the sample space $\Omega$ is the subset of nonsingular matrices in $\Mat_{n\times m}(\Q_p)$, the event space $\mathcal{F}$ is the $\sigma$-algebra generated by the following sets

$$\{(U'+p^d\Mat_{n\times n}(\Z_p))B(V'+p^d\Mat_{m\times m}(\Z_p))\mid d>0,U'\in\GL_n(\Z_p),V'\in\GL_m(\Z_p)\}$$
where $B\in\Mat_{n\times m}(\Q_p)$ is singular, and $G$ is the probability function such that 
for any $U_1',U_2'\in\GL_n(\Z_p), V_1',V_2'\in\GL_m(\Z_p)$ and $d>0$, the following two sets have the same probability:
$$\{(U_1'+p^d\Mat_{n\times n}(\Z_p))B(V_1'+p^d\Mat_{m\times m}(\Z_p))\},\{(U_2'+p^d\Mat_{n\times n}(\Z_p))B(V_2'+p^d\Mat_{m\times m}(\Z_p))\}.$$
Then $A$ has the form

$$A=U\diag_{n\times m}(p^{\lambda_1},\ldots,p^{\lambda_n})V$$
where the three random matrices on the right-hand side are independent, $U\in\GL_n(\Z_p), V\in\GL_m(\Z_p)$ are Haar distributed, and $\SN(A)=(\lambda_1,\ldots,\lambda_n)$ is randomly distributed over the set $\Sig_n$ of integer signatures of length $n$.
\end{cor}

\begin{example}\label{ex: nxm corner of invertible}
Let $1\le n\le m\le N$. Denote $A\in\Mat_{n\times m}(\Z_p)$ as the top $n\times m$ submatrix of a Haar-distributed element of $\GL_N(\Z_p)$. Then the distribution of $A$ is invariant under the action of $\GL_n(\Z_p)\times\GL_m(\Z_p)$. Also, when we take limit $N$ goes to infinity, we turn $A$ into the Haar matrix (i.e., the entries of $A$ are i.i.d Haar distributed over $\Z_p$). In this case, the distribution of $A$ is still invariant under the action of $\GL_n(\Z_p)\times\GL_m(\Z_p)$. See \Cref{thm: corner process and product process} for further discussions around this. 
\end{example}

The singular numbers of a given $p$-adic matrix are also related to its minors. Here by $k \times k$ minor, we mean any $k \times k$ matrix obtained by deleting rows and columns of the original matrix. The explicit form is as follows.

\begin{prop}\label{thm:submatrices_suffice}
Let $1 \leq n \leq m$ be integers and $A \in \Mat_{n \times m}(\Q_p)$ with $\SN(A) = (\lambda_1,\ldots,\lambda_n)$. Then for any $1 \leq k \leq n$,
\begin{equation}\label{eq:raleigh_submat}
\lambda_n+\ldots+\lambda_{n-k+1} = \min_{A' \text{ $k \times k$ minor of $A$}} \val_p(\det(A')).
\end{equation}
\end{prop}
\begin{proof}
See Proposition 2.4 of \cite{vanpeski2024padicdysonbrownianmotion}.
\end{proof}

We record a few other simple facts about singular numbers that will be useful.

\begin{prop}
\label{thm:multiply_smaller_dimension}
Let $n \leq m$, $A \in \Mat_{n \times m}(\Q_p)$, and $\kappa \in \Sig_m$. Then 
\[
|\SN(\diag_{n\times n}(p^{\kappa_1},\ldots,p^{\kappa_n}) A)| = |\SN(A)| + |\kappa|.
\]
\end{prop}
\begin{proof}
See Proposition 2.5 of \cite{vanpeski2024padicdysonbrownianmotion}.
\end{proof}

\begin{prop}\label{prop:singular numbers of sum}
Let $n \leq m$, $A=(a_{ij})_{1\le i\le n,1\le j\le m}\in\Mat_{n \times m}(\Q_p)$ with $\SN(A)=(\lambda_1,\ldots,\lambda_n)$. For any $1\le k\le n$ and $B=(b_{ij})_{1\le i\le n,1\le j\le m}\in\Mat_{n \times m}(\Q_p)$ with $\SN(B)_n>\SN(A)_k$, we have
$$\SN(A+B)_j=\SN(A)_j,\quad\forall k\le j\le n.$$
\end{prop}

\begin{proof}
By multiplication of $\GL_n(\Z_p)$ on the left side and $\GL_m(\Z_p)$ on the right side, without loss of generality, we may assume $A=\diag_{n\times m}(p^{\lambda_1},\ldots,p^{\lambda_n})$. Since all the entries of $B$ have values greater than $\lambda_k$, we can use the element $p^{\lambda_n}+b_{nn}$ to eliminate the last row and the $n$th column. This does not change the value of the entry on the intersection of the $(n-1)$th row and $(n-1)$th column, so we can repeat the operation to eliminate the $(n-1)$th row and the $(n-1)$th column, \ldots, until we get a diagonal with values $\lambda_k,\ldots,\lambda_n$. This ends the proof.
\end{proof}

\begin{prop}\label{thm:increase_decrease}
Let $n \leq m$, $A \in \Mat_{n \times m}(\Q_p)$, and suppose $B \in \Mat_m(\Q_p)$ has all singular numbers nonnegative (resp. nonpositive). Then $\SN(AB)_i \geq \SN(A)_i$ (resp. $\SN(AB)_i \leq \SN(A)_i$) for each $1 \leq i \leq n$. If $C \in \Mat_n(\Q_p)$ has nonnegative (resp. nonpositive) singular numbers, the same holds with $CA$ replacing $AB$. 
\end{prop}
\begin{proof}
By multiplying a power of $p$ to $A$, there is no loss we assume $A\in\Mat_{n\times m}(\Z_p)$. First, suppose $B \in \Mat_m(\Z_p)$ has all singular numbers nonnegative. Then, the identity map from $\Z_p^n$ to itself induces an abelian $p$-group embedding $\Z_p^n/A\Z_p^m\hookrightarrow\Z_p^n/AB\Z_p^m$. Therefore, the Young diagram of the latter includes the former, see \cite[Chapter 2, (4.3)]{Macdonald}. Hence by \cite[Remark 3]{vanpeski2024locallimitspadicrandom}, we identify $\SN(A)$ and $\SN(AB)$ as their cokernels and end the proof.

Next, suppose $C \in \Mat_n(\Z_p)$ has nonnegative singular numbers. Notice the map $\Z_p^n\rightarrow\Z_p^n,x\mapsto Cx$ induces an embedding $\Z_p^n/A\Z_p^m\hookrightarrow\Z_p^n/CA\Z_p^m$, we prove the proposition by the same method as above.

Finally, suppose $B \in \Mat_m(\Q_p)$ has all singular numbers nonpositive. Then by the above deduction, for each $1\le i\le n$, we have $\SN(A)_i=\SN(ABB^{-1})_i\ge\SN(AB)_i$; Suppose $C \in \Mat_n(\Q_p)$ has all singular numbers nonpositive. Then by the above deduction, for each $1\le i\le n$, we have $\SN(A)_i=\SN(C^{-1}CA)_i\ge\SN(CA)_i$.
\end{proof}

\begin{prop}\label{cor: interlacing}
Let $1\le i\le n-1$. Suppose $A^{(i)}\in\Mat_{(n-i+1)\times n}(\Q_p)$, and let $A^{(i+1)}\in\Mat_{(n-i)\times n}(\Q_p)$ be the submatrix of $A_i$ eliminating the first row. Then we have $\SN(A^{(i+1)})\prec_P\SN(A^{(i)})$, i.e., for all $1\le j\le n-i$,

$$\SN(A^{(i)})_{j+1}\le\SN(A^{(i+1)})_j\le\SN(A^{(i)})_j$$
\end{prop}

\begin{proof}
On one hand, by \Cref{thm:increase_decrease}, for sufficiently large $r$, we have
$$\SN(A^{(i+1)})_j=\SN(\diag_{(n-i+1)\times(n-i+1)}(p^r,1,\ldots,1)A^{(i)})_{j+1}\ge\SN(A^{(i)})_{j+1}$$
This is true because by \Cref{thm:submatrices_suffice}, when $r$ is sufficiently large, the minimal minor never involves the first row. 

On the other hand, by multiplying a power of $p$, there is no loss we assume $A^{(i)}\in\Mat_{(n-i+1)\times n}(\Z_p)$. By multiplying $\GL_{n-i}(\Z_p)$ on the left and $\GL_n(\Z_p)$ on the right, there is no loss we assume 
$$A^{(i)}=\begin{pmatrix}a & \beta \\ 0 & \diag_{(n-i)\times (n-1)}(p^{\mu_1},\ldots,p^{\mu_{n-i}})\end{pmatrix}$$
where $a\in\Z_p$, $\beta\in\Mat_{1\times(n-1)}(\Z_p)$, and $\mu_1\ge\ldots\ge\mu_{n-i}\ge 0$. By \Cref{thm:increase_decrease}, we have

\begin{align}\SN(A^{(i)})_j&=\SN(\diag_{n\times (n-i+1)}(1,p^{\mu_1},\ldots,p^{\mu_{n-i}})\begin{pmatrix}a & \beta \\ 0 & I_{n-1}\end{pmatrix})_j\\
&\ge\SN(\diag_{n\times(n-i+1)}(1,p^{\mu_1},\ldots,p^{\mu_{n-i}}))_j=\mu_j=\SN(A^{(i+1)})_j
\end{align}
which ends the proof.
\end{proof}

\begin{prop}\label{prop: weight of corners}
Let $A^{(i)},A^{(i+1)}$ be the same as \Cref{cor: interlacing}. Also, let $A^{(i+1)'}\in\Mat_{(n-i)\times n}(\Q_p)$ be the submatrix of $A_i$ eliminating the second row, and $A^{(i+2)}\in\Mat_{(n-i-1)\times n}(\Q_p)$ be the submatrix of $A_i$ eliminating the first and second row. Then we have 
\begin{equation}
|\SN(A^{(i)})|-|\SN(A^{(i+1)})|\ge|\SN(A^{(i+1)'})|-|\SN(A^{(i+2)})|.\end{equation}
\end{prop}
\begin{proof}
By multiplying $GL_{n-i+1}(\Z_p)$ on the left and $\GL_n(\Z_p)$ on the right, there is no loss we assume 
$$A^{(i)}=\begin{pmatrix}\alpha_1 & \beta_1 \\ \alpha_2 & \beta_2\\ 0 &\diag_{(n-i-1)\times (n-i-1)}(p^{\mu_1},\ldots,p^{\mu_{n-i-1}})  \end{pmatrix} 
$$
where $\alpha_1,\alpha_2\in\Mat_{1\times(i+1)}(\Q_p)$, and $\beta_1,\beta_2\in\Mat_{1\times(n-i-1)}(\Q_p)$. Then we have 
\begin{equation*}
|\SN(A^{(i)})|-|\SN(A^{(i+1)})|=\SN|\begin{pmatrix}\alpha_1\\ \alpha_2\end{pmatrix}|-|\SN(\alpha_2)|\ge|\SN (\alpha_1)|=|\SN(A^{(i+1)'})|-|\SN(A^{(i+2)})|.
\end{equation*}
\end{proof}

\begin{rmk}
if $i=n-1$, we delete the form $|\SN(A^{(i+2)})|$ in the above inequality. From now on, we no longer restate this elimination since it is clear from the text.
\end{rmk}

\begin{rmk}
Throughout the paper, we require the measures over each matrix to be left- and right-invariant under the multiplication of $\GL_n(\Z_p)$. However, the same method still works when the measures are only left-invariant or only right-invariant. The example of the product of two random matrices $A_1A_2$ whose measures are invariant under left-multiplication of $\GL_n(\Z_p)$ will explain everything. We write decomposition

$$A_1=U_1B_1,\quad A_2=U_2B_2$$
where $U_1,U_2,B_1,B_2$ are independent random matrices, and $U_1,U_2\in\GL_n(\Z_p)$ are Haar distributed. We write $n$-dimensional random variables $\SN(A_1)=(\mu_1(1),\ldots,\mu_n(1))$, and $\SN(A_2)=(\mu_1(2),\ldots,\mu_n(2))$. In this case, we have ($\stackrel{d}{=}$ here means ``have the same distribution")

\begin{align}
\SN(A_1A_2)&=\SN(U_1B_1U_2B_2)=\SN(B_1U_2B_2)\nonumber\\
&\stackrel{d}{=}\SN(\diag_{n\times n}(p^{\mu_1(1)},\ldots,p^{\mu_n(1)})U_2\diag_{n\times n}(p^{\mu_1(2)},\ldots,p^{\mu_n(2)}))\nonumber\\
&\stackrel{d}{=}\SN(U_1\diag_{n\times n}(p^{\mu_1(1)},\ldots,p^{\mu_n(1)})V_1U_2\diag_{n\times n}(p^{\mu_1(2)},\ldots,p^{\mu_n(2)})V_2)
\end{align}
where the $V_1,V_2\in\GL_n(\Z_p)$ are Haar distributed, independent with all the other matrices. The second row holds because for all fixed $B_1,B_1',B_2,B_2'$ such that $\SN(B_1)=\SN(B_1'),\SN(B_2)=\SN(B_2')$, we have $\SN(B_1U_2B_2)\stackrel{d}{=}\SN(B_1'U_2B_2')$. Hence, we are back to the familiar situation: the measures are invariant on both sides.
\end{rmk}


\section{Rigorous description of the heuristic and proof of i.i.d case}\label{Proof}

In this section, we provide our method based on the backgrounds in the last section. Suppose that $A_1, A_2,\ldots\in\GL_n(\Q_p)$ are independent random matrices, and their measures are invariant under the left- and right-multiplication of $\GL_n(\Z_p)$.  For all $k\ge 1, \mu\in\Sig_n$, denote $\mathbf{P}_k(\mu):=\mathbf{P}(\SN(A_k)=\mu)$. In accordance with \Cref{cor: invariant measure}, let $(\Omega_k,\mathcal{F}_k,G_k)$ denote the probability space of $A_k$.

\begin{defi}\label{defi: probability space}
Define the product space $(\Omega,\mathcal{F},G)$,
corresponding to the sequence $(A_1,A_2,\ldots)$, as:
$$\Omega=\Omega_1\times\Omega_2\times\ldots.$$
This construction of probability space is achievable via Kolmogrov's extension theorem.
\end{defi}

For every $1\le i\le n$, let $A_k^{(i)}\in\Mat_{(n-i+1)\times n}(\Q_p)$ denote the submatrix of $A_k$, which is the last $n-i+1$ columns and all $n$ rows. Then by \Cref{cor: interlacing}, we have the interlacing relation:
$$\SN(A_k^{(n)})\prec_P\ldots\prec_P\SN(A_k^{(1)}).$$

Following the approach outlined in the Introduction section, we construct a random sequence of integer signatures $\lambda(0)=(0[n]),\lambda(1),\lambda(2),\ldots\in\Sig_n$ over $\Omega$ by induction. First, $\lambda(1)=\SN(A_1)$ is the singular numbers of $A_1$; If we already know $\lambda(k-1)$, then $\lambda(k)=(\lambda_1(k),\ldots,\lambda_n(k))$ would be the singular numbers of the matrix
$$\diag_{n\times n}(p^{\lambda_1(k-1)},\ldots,p^{\lambda_n(k-1)})A_{k-1}.$$

Since the measure is invariant under the action of $\GL_n(\Z_p)$ on both sides, we know the joint distribution of $(\lambda(1),\lambda(2),\ldots)$ is the same as the joint distribution of  $(\SN(A_1),\SN(A_1A_2),\ldots)$. On the other hand, let $v(0)=(0[n]),v(1), v(2),\ldots$ be the random sequence in $\Z^n$ such that for all $k\ge 0$ and $1\le j\le n$,
\begin{equation}\label{eq: definition of v}
v_i(k)+\ldots+v_n(k)=\sum_{j=1}^k |\SN(A_j^{(i)})|
\end{equation}
where for every $1\le j\le n$,  $A_i^{(j)}$ is the submatrice of $A_i$, consisting of the last $n-j+1$ rows and all $n$ columns. In other words, $v$ is generated by adding singular numbers of the corners, and no matrix product is involved throughout the calculation. In this case, the relation $v_1(k)\ge\ldots\ge v_n(k)$ might not hold, and $v(k)=(v_1(k),\ldots,v_n(k))$ does not need to be in $\Sig_n$.

We now introduce a key definition, as noted in the Introduction section.

\begin{defi}\label{defi: split}

Suppose $\mu=(\mu(0),\mu(1),\mu(2),\ldots)$ is random over the probability space $\Omega=\Omega_1\times\Omega_2\times\ldots$ defined in \Cref{defi: probability space}, where $\mu(k)\in\Z^n$ for all $k\ge 0$ . We say the sequence $\mu$ is \textbf{split}, if for all $1\le i\le n-1$,
$$\lim_{k\rightarrow\infty}\mu_i(k)-\mu_{i+1}(k+1)+\SN(A_{k+1})_n=+\infty.$$
\end{defi}

\begin{example}\label{example: split} 
For many cases, including the examples below, we would see that the sequence $v$ we defined in \eqref{eq: definition of v} is split almost surely.

\begin{enumerate}[left=0pt]
\item\label{item: i.i.d split} Suppose $A_1,A_2,\ldots$ are i.i.d. random matrices. In this case, the random variable $v(k)$ is the sum of i.i.d random variables of dimension $n$. If we have $\sum_{M\in\Z}\mathbf{P}(\SN(A_1)=M[n])=1$, it is clear that $v$ is never split. Otherwise, if there exists $\mu^{(1)}=(\mu_1^{(1)},\ldots\mu_n^{(1)})\in\Sig_n$ such that $\mathbf{P}(\SN(A_1)=\mu^{(1)})>0$ and $\mu_1^{(1)}>\mu_n^{(1)}$, we would have

\begin{equation}\label{eq: inequality of corners}
\E|\SN(A_1^{(n)})|<\E|\SN(A_1^{(n-1)})|-\E|\SN(A_1^{(n)})|<\ldots<\E|\SN(A_1^{(1)})|-\E|\SN(A_1^{(2)})|\end{equation}
as long as these expectations exist. (We will return and prove this inequality at the end of \Cref{Proof}, where we continue our proof of \Cref{cor: i.i.d.}.) Thus, we apply the strong law of large numbers and prove
$$\frac{v_i(k)}{k}\rightarrow \E|\SN(A_1^{(i)})|-\E|\SN(A_1^{(i+1)})|\quad\forall 1\le i\le n-1,\quad\frac{v_n(k)}{k}\rightarrow \E|\SN(A_1^{(n)})|$$
almost surely. Thus, the sequence $v$ is split almost surely.

\item \label{item: Roger setting}Under the settings of Van Peski \cite[Theorem 1.1]{Roger}, for all $k\ge 1$, $A_k$ is the top $n\times n$ matrix of a Haar-distributed element of $\GL_{N_k}(\Z_p)$, where $N_1,N_2,\ldots\in\Z\cup\{\infty\}$ with $N_k>n$ for all $k\ge 1$. In this case, the sequence $v$ we defined coincides with the previous paper \cite[Definition 15]{Roger}, hence the assertion that the sequence $v$ is almost surely split has already been proved. See \Cref{rmk: Roger setting split} for details.
\end{enumerate}
\end{example}

The following milestone responds to the heuristic in the introduction section, which proves that the split property implies the difference between $\lambda(k)$ and $v(k)$ is negligible when studying asymptotics.

\begin{thm}\label{thm: main theorem}

Suppose that $v$ is split. Then, the supremum of the difference
$$\sup_{k\in\Z_{\ge 0}}|\lambda_i(k)-v_i(k)|$$
is bounded for every $1\le i\le n$.

\end{thm}

As preparation for our proof, we construct random sequences $\lambda^{(1)}(k)=v(k),\ldots,\lambda^{(n)}(k)=\lambda(k)\in\Z_{\ge 0}^n$, all defined over $\Omega$. The construction starts from $\lambda^{(1)}$, and ends at $\lambda^{(n)}$ one by one. For all $1\le j\le n$, we construct $\lambda^{(j)}(k)$ by induction over $k$. When $k=0$, we let $\lambda^{(j)}(0)=(0,0,\ldots,0)$; If we already know $\lambda^{(j)}(k)$, then:
\begin{enumerate}[left=0pt]
\item $(\lambda_{n-j+1}^{(j)}(k+1),\ldots,\lambda_n^{(j)}(k+1))\in\Sig_j=\SN(\diag_{j\times j}(p^{\lambda_{n-j+1}^{(j)}(k)},\ldots,p^{\lambda_n^{(j)}(k)})A^{(n-j+1)}_k)$;

\item For all $1\le i\le n-j$, set $\lambda_i^{(j)}(k+1)=\lambda_i^{(j-1)}(k+1)=\ldots=\lambda_i^{(1)}(k+1)=v_i(k+1)$.
\end{enumerate}
Under the above construction, we must have $\lambda_{n-j+1}^{(j)}(k)\ge\ldots\ge\lambda_n^{(j)}(k)$ for all $k$, but there is no such descending relation for $\lambda_1^{(j)}(k+1),\ldots,\lambda_{n-j}^{(j)}(k+1)$. We will prove by induction on $j$ that when $v$ is split,

$$\sup_{k\in\Z_{\ge 0}}|\lambda_i^{(j)}(k)-v_i(k)|$$
is finite for all $i$. When $j=n$, we are done.

\quad

The sequences $\lambda^{(j)}$ and $\lambda^{(j+1)}$ satisfy the following inequality:

\begin{prop}\label{prop: neighbour relation}
For all $k\ge 0$ and $1\le i\le j$, we have
$$\lambda_{n-j}^{(j+1)}(k)\ge\lambda_{n-j}^{(j)}(k),\quad\lambda_{n-j+i}^{(j+1)}(k)\le\lambda_{n-j+i}^{(j)}(k).$$
\end{prop}

\begin{proof}
We prove this by induction over $k$. The base case $k=0$ follows since $\lambda^{(l)}(0)=0[n]$. Suppose that the above holds for some $k$. Then we have for all $1\le i\le j$,

\begin{align}
\lambda_{n-j+i}^{(j+1)}(k+1)&=\SN(\diag_{(j+1)\times(j+1)}(p^{\lambda_{n-j}^{(j+1)}(k)},\ldots,p^{\lambda_n^{(j+1)}(k)})A^{(n-j)}_k)_{i+1}\nonumber\\
&\le\SN(\diag_{j\times j}(p^{\lambda_{n-j+1}^{(j+1)}(k)},\ldots,p^{\lambda_n^{(j+1)}(k)})A^{(n-j+1)}_k)_i\nonumber\\
&\le\SN(\diag_{j\times j}(p^{\lambda_{n-j+1}^{(j)}(k)},\ldots,p^{\lambda_n^{(j)}(k)})A^{(n-j+1)}_k)_i\nonumber\\
&=\lambda_{n-j+i}^{(j)}(k+1).
\end{align}
The inequality on the second row is due to \Cref{cor: interlacing}, and the inequality on the third row is due to the induction hypothesis. Also, notice that $|\lambda^{(j+1)}(k+1)|=|\lambda^{(j)}(k+1)|=\sum_{i=1}^{k+1}|\SN(A_i)|$, and $\lambda_i^{(j+1)}(k+1)=\lambda_i^{(j)}(k+1)=v_i(k+1)$ for all $1\le i\le n-j-1$. Thus we have

\begin{equation}\label{eq: equal sum}
\lambda_{n-j}^{(j+1)}(k+1)+\sum_{1\le i\le j}\lambda_{n-j+i}^{(j+1)}(k+1)=\lambda_{n-j}^{(j)}(k+1)+\sum_{1\le i\le j}\lambda_{n-j+i}^{(j)}(k+1).
\end{equation}
This imply $\lambda_{n-j}^{(j+1)}(k+1)\ge\lambda_{n-j}^{(j)}(k+1)$. Therefore, the proposition is also true for $k+1$.
\end{proof}

The following lemma connects the increment of singular numbers with matrix corners.

\begin{lemma}\label{lem: equal difference}
Suppose $1\le j\le n-1$, $(\lambda_{n-j},\ldots,\lambda_n)\in\Sig_{j+1}$, $A^{(n-j)}\in\GL_{(j+1)\times n}(\Q_p)$ be full rank, and $A^{(n-j+1)}\in\GL_{j\times n}(\Q_p)$ be the submatrix of $A^{(n-j)}$, consisting of the last $j$ rows and all $n$ columns. Suppose we have the inequality

\begin{equation}\label{eq: first and second}
\SN(\diag_{(j+1)\times(j+1)}(p^{\lambda_{n-j}},\ldots,p^{\lambda_n})A^{(n-j)})_2<\lambda_{n-j}+\SN(A^{(n-j)})_{j+1}.
\end{equation}
Then we have 
$$\SN(\diag_{(j+1)\times(j+1)}(p^{\lambda_{n-j}},\ldots,p^{\lambda_n})A^{(n-j)})_1-\lambda_{n-j}=|\SN(A^{(n-j)})|-|\SN(A^{(n-j+1)})|.$$
\end{lemma}

\begin{proof}
By multiplying a power of $p$ to $A^{(n-j)}$, there is no loss we assume $\SN(A^{(n-j)})_{j+1}=0$, and $A^{(n-j)}\in\Mat_{(j+1)\times n}(\Z_p)$ has nonnegative singular numbers. In this case, we have $\SN(\diag_{(j+1)\times(j+1)}(p^{\lambda_{n-j}},\ldots,p^{\lambda_n})A^{(n-j)})_2<\lambda_{n-j}$. Let 
$$B^{(n-j)}:=\diag_{(j+1)\times(j+1)}(p^{\lambda_{n-j}},\ldots,p^{\lambda_n})A^{(n-j)}, B^{(n-j+1)}:=\diag_{j\times j}(p^{\lambda_{n-j+1}},\ldots,p^{\lambda_n})A^{(n-j+1)}$$
so $B^{(n-j+1)}$ becomes the submatrix of $B^{(n-j)}$, moving the first row. First, by \Cref{prop:singular numbers of sum}, we have
$$\lambda_{n-j}>\SN(B^{(n-j)})_2=\SN(\begin{pmatrix}0_{1\times n}\\ B^{(n-j)}\end{pmatrix})_2=\SN(B^{(n-j+1)})_1.$$

Next, for any $j\times j$ minor $B'$ of $B^{(n-j)}$ such that the first row is picked, we have

\begin{align}\val_p(\det(B'))&\ge\lambda_{n-j}+\inf_{A' \text{ $(j-1) \times (j-1)$ minor of $B^{(n-j+1)}$}}\val_p(\det(A'))\nonumber\\
&>\SN(B^{(n-j+1)})_1+\inf_{A' \text{ $(j-1) \times (j-1)$ minor of $B^{(n-j+1)}$}}\val_p(\det(A'))\nonumber\\
&=\inf_{A' \text{ $j \times j$ minor of $B^{(n-j+1)}$}}\val_p(\det(A'))=\SN(B^{(n-j+1)}).\nonumber\\
\end{align}
Hence the minimal value of $j\times j$ minor of $B^{(n-j)}$ is $\SN(B^{(n-j+1)})$, and we have
\begin{align}
\SN(B^{(n-j)})_1&=|\SN(B^{(n-j)})|-|\SN(B^{(n-j+1)})|\nonumber\\
&=(\lambda_{n-j}+\ldots+\lambda_n)+|\SN(A^{(n-j)})|-(\lambda_{n-j+1}+\ldots+\lambda_n)-|\SN(A^{(n-j+1)})|\nonumber\\
&=\lambda_{n-j}+|\SN(A^{(n-j)})|-|\SN(A^{(n-j+1)})|.\nonumber\\
\end{align}
This ends the proof of the lemma.
\end{proof}

\begin{proof}[Proof of \Cref{thm: main theorem}]
Based on \Cref{lem: equal difference}, we show by induction over $j$ that when $v$ is split,  the difference between $\lambda^{(j)}$ and $v$ is bounded, i.e.,

$$\sup_{k\ge 1}\sup_{1\le i\le  j}|\lambda_{n-i+1}^{(j)}(k)-v_{n-i+1}(k)|<\infty.$$

The base case $j=1$ follows since $\lambda^{(1)}=v$; Now, suppose the above holds for some $j$. Since the difference between $\lambda^{(j)}$ and $v$ is bounded, $\lambda^{(j)}$ is also split. On one hand, notice that for $k$ such that $\lambda_{n-j+1}^{(j+1)}(k+1)-\SN(A_{k+1})_n<\lambda_{n-j}^{(j+1)}(k)$,

$$\lambda_{n-j}^{(j+1)}(k+1)-\lambda_{n-j}^{(j+1)}(k)=|\SN(A_k^{(n-j)})|-|\SN(A_k^{(n-j+1)})|=v_{n-j}(k+1)-v_{n-j}(k).$$

On the other hand, since $\lambda^{(j)}$ is split, we have $\lambda_{n-j+1}^{(j)}(k+1)-\SN(A_{k+1})_n<\lambda_{n-j}^{(j)}(k)$ holds except for finitely many $k$. Therefore, by \Cref{prop: neighbour relation}, the inequality 

$$\lambda_{n-j+1}^{(j+1)}(k+1)-\SN(A_{k+1})_n\le\lambda_{n-j+1}^{(j)}(k+1)-\SN(A_{k+1})_n<\lambda_{n-j}^{(j)}(k)\le\lambda_{n-j}^{(j+1)}(k)$$ also holds except for finitely many $k$. Hence, the increments $\lambda_{n-j}^{(j+1)}(k+1)-\lambda_{n-j}^{(j+1)}(k)$ and $v_{n-j}(k+1)-v_{n-j}(k)$ are equal except for finitely many $k$, and

$$\sup_{k\ge 1}|\lambda_{n-j}^{(j+1)}(k)-v_{n-j}(k)|<\infty,\quad\sup_{k\ge 1}(\lambda_{n-j}^{(j+1)}(k)-\lambda_{n-j}^{(j)}(k))<\infty.$$

Also, by \Cref{prop: neighbour relation} and \eqref{eq: equal sum}, we have for all $1\le i\le j$,

$$\sup_{k\ge 1}(\lambda_{n-i+j}^{(j)}(k)-\lambda_{n-j+i}^{(j+1)}(k))\le\sup_{k\ge 1}(\lambda_{n-j}^{(j+1)}(k)-\lambda_{n-j}^{(j)}(k))<\infty.$$
This implies that the difference between $\lambda^{(j+1)}$ and $\lambda^{(j)}$ is bounded, and by the induction hypothesis, the theorem is also true for $j+1$, which ends our proof.
\end{proof}

\begin{example}
If $v$ is not split, our conclusion that the difference between $\lambda$ and $v$ is bounded might not hold. Consider the case $n=2$, and for all $k\ge 1$, we set $A_k=\begin{pmatrix}1 & 0\\0 & p^{2^{k-1}} \end{pmatrix}$. Then, we can verify that
$$\lambda_1(k)=2^{k-1}+2^{k-3}+\ldots+2^{k-1-2[\frac{k-1}{2}]},\quad\lambda_2(k)=2^{k-2}+2^{k-4}+\ldots+2^{k-2[\frac{k}{2}]}$$
$$v_1(k)=0,\quad v_2(k)=2^{k-1}+2^{k-2}+\ldots+1=2^k-1.$$
Therefore, this gives a counterexample.
\end{example}

\begin{rmk}\label{rmk: method}
\Cref{thm: main theorem} confirms the heuristic in the introduction section that for most cases, the asymptotic of independent matrix products can be derived from simply adding the corners, which turns the problem into studying the sum of $n$-dimensional random variables, and covers a large range of distributions including the i.i.d. case. As a summary, to study the asymptotic of the joint distribution of $\SN(A_1),\SN(A_1A_2),\ldots$, we go through the following steps:

\begin{enumerate}[left=0pt]
\item For each $k\ge 1$, we find the joint distribution of the weights $|\SN(A_k^{(1)})|,\ldots,|\SN(A_k^{(n)})|$. It is worth mentioning that \Cref{cor: joint distribution of corner} in the appendix provides the joint distribution of $\SN(A_k^{(1)}),\ldots,\SN(A_k^{(n)})$ in the form of Hall-Littlewood polynomials;
\item We add those independent random variables and get the sequence $v=(v(0),v(1),\ldots)$. This step concerns the sum of independent random variables and has nothing to do with random matrix products. In many cases, the strong law of large numbers can be applied;
\item As long as the sequence $v$ is split, the joint distribution 

$$\lambda(0)=0[n],\lambda(1)=\SN(A_1),\lambda(2)=\SN(A_1A_2),\ldots$$
shares the same asymptotic as $v$.
\end{enumerate}
\end{rmk}

\begin{proof}[Proof of \Cref{cor: i.i.d.}]\label{Proof of Theorem 1.1}
We now continue the proof of the i.i.d case in \eqref{item: i.i.d split} of \Cref{example: split}. Suppose there exists $\mu^{(1)}\in\Sig_n$ such that $\mathbf{P}_1(\mu^{(1)}):=\mathbf{P}(\SN(A_1)=\mu^{(1)})>0$ and $\mu_1^{(1)}>\mu_n^{(1)}$. In this case, if we assume the inequality \eqref{eq: inequality of corners} holds, then the sequence $v$ is split almost surely, and \Cref{thm: main theorem} implies the asymptotic of the sequence $\lambda$ is almost surely the same as $v$. The sum of i.i.d. random variables would imply the following:

\begin{enumerate}[left=0pt]
\item (Strong law of large numbers) Suppose that for all $1\le i\le n$, the expectation $\E|\SN(A_1^{(i)})|<\infty$ exists. Then for all $1\le i\le n$, we have
$$\frac{v_i(k)+\ldots+v_n(k)}{k}\rightarrow \E|\SN(A_1^{(i)})|$$
almost surely as $k$ goes to infinity.

\item(Central limit theorem) In addition to the conditions above, assume that for all $1\le i\le n$, the expectation $\E|\SN(A_1^{(i)})^2|<\infty$ exists. Then we have weak convergence
$$\frac{(v_1(k)+\ldots+v_n(k)-k\E|\SN(A_1^{(1)})|,\ldots,v_n(k)-k\E|\SN(A_1^{(n)})|)}{\sqrt{k}}\Rightarrow N(0,\Sigma)$$
as $k$ goes to infinity.
\end{enumerate}
Replace $v$ by $\lambda$, and a linear transform for the multivariable normal distribution yields the result. Therefore, the only step left for us is to prove the inequality \eqref{eq: inequality of corners}. By \Cref{prop: weight of corners}, we have for all $1\le i\le n-1$,
\begin{equation}\label{eq: i.i.d_inequality of weight}
|\SN(A_1^{(i)})|-|\SN(A_1^{(i+1)})|\ge|\SN(A_1^{(i+1)'})|-|\SN(A_1^{(i+2)})|\end{equation}
where $A_1^{(i+1)'}$ denotes the submatrix of $A$ consisting of the $i^{\text{th}},(i+2)^{\text{th}},(i+3)^{\text{th}}\ldots,n^{\text{th}}$ rows and all $n$ columns. $A_1^{(i+1)'}$ has the same distribution as $A_1^{(i+1)}$, because the two matrices $A_1$ and
$$\begin{pmatrix}I_{i-1} &  &  & \\  & 0 & 1 &  \\  & 1 & 0 & \\ &  &  & I_{n-i-1}\end{pmatrix}A_1$$
share the same distribution, and $A^{(i+1)'}$ is the submatrix of the latter consisting of the last $n-i$ rows and $n$ columns. Set $M=1+p^D\Mat_{n\times n}(\Z_p)$, where $D>\mu_1^{(1)}-\mu_n^{(1)}$ is sufficiently large. By \Cref{prop:singular numbers of sum}, the matrices in the set
$$\setlength{\arraycolsep}{2pt}
S=M\begin{pmatrix}\diag_{(i-1)\times(i-1)}(p^{\mu_2^{(1)}},\ldots,p^{\mu_i^{(1)}}) &  &  & \\
 & p^{\mu_n^{(1)}} & 0 &  \\
 & p^{\mu_n^{(1)}} & p^{\mu_1^{(1)}} & \\
 &  &  & \diag_{(n-i-1)\times(n-i-1)}(p^{\mu_{i+1}^{(1)}},\ldots,p^{\mu_{n-1}^{(1)}})
\end{pmatrix}M$$
all have singular numbers $\mu^{(1)}=(\mu_1^{(1)},\ldots,\mu_n^{(1)})$. The probability $\mathbf{P}(A_1\in S)>0$ is nonzero, and whenever this happens, \eqref{eq: i.i.d_inequality of weight} holds strict inequality. 
Therefore, by taking expectations on both sides of \eqref{eq: i.i.d_inequality of weight}, we have $\E|\SN(A_1^{(i)})|-\E|\SN(A_1^{(i+1)})|>\E|\SN(A_1^{(i+1)})|-\E|\SN(A_1^{(i+2)})|$ for all $1\le i\le n-1$.
\end{proof}

\section{Random matrix product over split reductive group}\label{reductive group}

Let $\G$ be a split reductive group over $\Q_p$. In this section, we provide an outlook on how our heuristic extends to random matrices over $G:=\G(\Q_p)$ with distribution invariant under the left- and right-multiplication of $K:=\G(\Z_p)$, as discussed in \Cref{rmk: reductive group}. In particular, for the case $\G=\GSp_{2n}$, we prove the strong law of large numbers and the central limit theorem, which turns out to be the splitting image of the $\GL_n$ case. For representation theory backgrounds in reductive groups, Hecke algebra, and root systems, see \cite{bruhat1972groupes}\cite{bruhat1984groupes}. 

Given $n\ge 1$ and $A\in G=\G(\Q_p)$, by Cartan decomposition there exists $U,V\in K$ such that $A=Ud(\lambda)V$ where $d(\lambda)$ is the normal form associated with the cocharacter $\lambda\in P^+$. We still denote $\SN(A)=\lambda$ as the \textbf{singular numbers} of $A$. In the typical case $\G=\GL_n$ that has been previously discussed, we have $P^+=\Sig_n$, and $d(\lambda)$ is the diagonal matrix $\diag_{n\times n}(p^{\lambda_1},\ldots,p^{\lambda_n})$, where the integers $\lambda_1\ge\ldots\ge\lambda_n$ are the parts of $\lambda$.
\begin{example}\label{ex: GSp}
Let $\G=\GSp_{2n}$ be the general symplectic group, $J_n=\renewcommand{\arraystretch}{0.5}\setlength{\arraycolsep}{2pt}\begin{pmatrix}& & 1\\ & \iddots & \\ 1 & &\end{pmatrix}\in\Mat_{n\times n}(\Z_p)$, and
$$G=\GSp_{2n}(\Q_p)=\{A\in\GL_{2n}(\Q_p)\mid A\begin{pmatrix} & J_n \\ -J_n & \end{pmatrix}A^T=\begin{pmatrix}  & \mu(A)J_n \\ -\mu(A)J_n & \end{pmatrix},\mu(A)\in\Q_p^\times\},$$
$$K=\GSp_{2n}(\Z_p)=G\cap\GL_{2n}(\Z_p).$$
Then for every $A\in G$, there exists $U,V\in K$ such that
$$A=U\diag_{2n\times 2n}(p^{\lambda_1},\ldots,p^{\lambda_{2n}})V$$
where $(\lambda_1\ge\ldots\ge\lambda_{2n})\in\Sig_{2n}$, and $\lambda_1+\lambda_{2n}=\lambda_2+\lambda_{2n-1}=\ldots=\lambda_n+\lambda_{n+1}$. We refer the integers $\lambda_i$ as the \textbf{singular numbers} of $A$, denoted as $\SN(A)=(\lambda_1,\ldots,\lambda_{2n})$.
\end{example}
Since $K$ is compact, we can define the multiplicative Haar measure over $K$. Suppose $A\in G$ is random with distribution invariant under left- and right-multiplication of $K$. Analogous to \Cref{cor: invariant measure}, $A$ has the form
$$A=Ud(\lambda)V$$
where the three random matrices on the right-hand side are independent, $U,V\in K$ are Haar distributed, and $\lambda$ is randomly distributed within the root system $P^+$.

Similar to the case $\G=\GL_n$ we discussed in the previous sections, we now investigate the asymptotic behavior of the sequence $\SN(A_1),\SN(A_1A_2),\SN(A_1A_2A_3),\ldots$, where $A_1,A_2,\ldots\in G$ are independent matrices with distribution invariant under left- and right-multiplication of $K$. Given $A_1,\ldots,A_k$, we write $\lambda(k)=\SN(A_1\cdots A_k)$, and we aim to determine the distribution $\SN(A_1\cdots A_{k+1})$. Taking advantage of the fact that the distribution is invariant under left- and right-multiplication by $K$, we are allowed to normalize the matrix $A_1\cdots A_k$ without changing the distribution and thus examine the singular numbers
$$\lambda(k+1)=\SN(d(\lambda(k))A_{k+1}).$$
when $k$ goes to infinity, it is reasonable to expect $\lambda(k)$ exhibit ``split" behaviors similar to the case $\G=\GL_n$. Therefore, multiplying a new random matrix $A_{k+1}$ could have a straightforward interpretation, which simplifies the problem into studying the independent sum of random variables. For instance, consider the case $\G=\GSp_{2n},G=\GSp_{2n}(\Q_p),K=\GSp_{2n}(\Z_p)$ from \Cref{ex: GSp}. When $k$ goes to infinity, one again looks forward to seeing the distance between the parts of $\lambda(k)=(\lambda_1(k),\ldots,\lambda_{2n}(k))$ also goes to infinity, i.e., $\lambda_1(k)$ is much larger than $\lambda_2(k)$, $\lambda_2(k)$ is much larger than $\lambda_3(k)$, and so on. Hence, by \Cref{thm:submatrices_suffice} the two sequences below share the same asymptotic:

\begin{enumerate}[left=0pt]
\item \label{item: matrix product_GSp}The sequence $\lambda(1),\lambda(2),\ldots$, where the $2n$-dimension joint variable $\lambda(k)$ is the singular numbers of matrix product $A_1\cdots A_k$;
\item \label{item: add corner_GSp} The sequence $v(1),v(2),\ldots$, where the $2n$-dimension joint variable $v(k)=(v_1(k),\ldots,v_{2n}(k))\in\Z^{2n}$ satisfies 
$$v_i(k)+\ldots+v_n(k)=\sum_{j=1}^k |\SN(A_j^{(i)})|.$$
\end{enumerate}  
Here $A_1^{(i)}\in\Mat_{(2n-i+1)\times n}(\Q_p)$ denote the submatrix of $A_1$ consisting of the last $2n-i+1$ rows and all $2n$ columns, and we regard $\SN(A_1^{(i)})$ as the entries of the diagonal form under the action of $\GL_{2n-i+1}(\Z_p)\times \GL_{2n}(\Z_p)$, which comes from the typical $\G=\GL_{2n}$ case. Define the product space $(\Omega,\mathcal{F},G)$,
corresponding to the sequence $(A_1,A_2,\ldots)$ as
$\Omega=\Omega_1\times\Omega_2\times\ldots$ via Kolmogrov's extension theorem. As an analogy of \Cref{thm: main theorem}, we have the following:

\begin{thm}
Suppose that $v$ is split, i.e., for all $1\le i\le 2n-1$, 
$$\lim_{k\rightarrow\infty}v_i(k)-v_{i+1}(k+1)+\SN(A_{k+1})_{2n}=+\infty.$$Then, the supremum of the difference $\sup_{k\in\Z_{\ge 0}}|\lambda_i(k)-v_i(k)|$
is bounded for every $1\le i\le 2n$.
\end{thm}
The proof of the above is similar to \Cref{thm: main theorem}, and we omit it here. In particular, when the matrices $A_1,A_2,\ldots$ are i.i.d, we prove the following theorem: 

\begin{thm}\label{thm: i.i.d for GSp}Let $A_1,A_2,\ldots\in\GSp_{2n}(\Q_p)$ be i.i.d. random matrices whose distributions are invariant under both left- and right-multiplication by $\GSp_{2n}(\Z_p)$. Let $\lambda(k)=(\lambda_1(k),\ldots,\lambda_{2n}(k))=\SN(A_1\ldots A_k)$. For every $1\le i\le 2n$, let $A_1^{(i)}\in\Mat_{(2n-i+1)\times (2n)}(\Q_p)$ denote the submatrix of $A_1$ consisting of the last $2n-i+1$ rows and all $2n$ columns. 
\begin{enumerate}[left=0pt]
\item (Strong law of large numbers) Suppose that for all $1\le i\le 2n$, the expectation $\E|\SN(A_1^{(i)})|<\infty$ exists. Then we have
\begin{align}(\frac{\l_1(k)}{k},\ldots,\frac{\l_{2n}(k)}{k})\rightarrow &(\E|\SN(A_1^{(1)})|-\E|\SN(A_1^{(2)})|,\ldots,\nonumber\\
&\E|\SN(A_1^{(2n-1)})|-\E|\SN(A_1^{(2n)})|,\E|\SN(A_1^{(2n)})|)\nonumber\end{align}
almost surely as $k$ goes to infinity. 

\item(Central limit theorem) In addition to the conditions above, assume that for all $1\le i\le 2n$, the expectation $\E|\SN(A_1^{(i)})^2|<\infty$ exists. Then we have weak convergence
$$\frac{(\l_1(k)-k\E|\SN(A_1^{(1)})|+k\E|\SN(A_1^{(2)})|,\ldots,\l_{2n}(k)-k\E|\SN(A_1^{(2n)})|)}{\sqrt{k}}\Rightarrow N(0,L\Sigma L^T)$$
as $k$ goes to infinity, where $L^T\in\Mat_{2n\times 2n}(\Z_p)$ is the transpose of the matrix
$$L=\begin{pmatrix}1 & -1 & & &\\ & 1 & -1 & &\\ & & \ddots & \ddots &\\ & & & 1 & -1\\ & & & & 1\end{pmatrix}$$
and $\Sigma=\Cov_{1\le i,j\le 2n}(|\SN(A_1^{(i)})|,|\SN(A_1^{(j)})|)$ is the covariance matrix.
\end{enumerate}
\end{thm}
\begin{proof}Our approach is very similar to the $\GL_n$ case. Assuming there exists $\mu^{(1)}=(\mu_1^{(1)},\ldots,\mu_{2n}^{(1)})$ such that 
\begin{equation}(\mu_1^{(1)}\ge\ldots\ge\mu_{2n}^{(1)})\in\Sig_{2n},\mu_1^{(1)}+\mu_{2n}^{(1)}=\ldots=\mu_n^{(1)}+\mu_{n+1}^{(1)},\mu_1^{(1)}>\mu_{2n}^{(1)},\mathbf{P}(\SN(A_1)=\mu^{(1)})>0,\end{equation}
we only need to prove for all $1\le i\le 2n-1$, the strict inequality 
\begin{equation}
\E|\SN(A_1^{(i)})|-\E|\SN(A_1^{(i+1)})|>\E|\SN(A_1^{(i+1)})|-\E|\SN(A_1^{(i+2)})|
\end{equation}
must hold. By \Cref{prop: weight of corners}, we have the following inequality
\begin{equation}\label{eq: i.i.d_inequality of weight for GSp}
|\SN(A_1^{(i)})|-|\SN(A_1^{(i+1)})|\ge|\SN(A_1^{(i+1)'})|-|\SN(A_1^{(i+2)})|\end{equation}
where $A_1^{(i+1)'}$ denotes the submatrix of $A$ consisting of the $i^{\text{th}},(i+2)^{\text{th}},(i+3)^{\text{th}},\ldots,(2n)^{\text{th}}$ rows and all $2n$ columns. $\SN(A_1^{(i+1)'})$ has the same distribution as $\SN(A_1^{(i+1)})$, because the two matrices $A_1$ and
$$\begin{cases}
\diag_{(2n)\times(2n)}(I_{i-1},\begin{pmatrix} 0 & 1\\ 1 & 0\end{pmatrix},I_{2n-2i-2},\begin{pmatrix} 0 & 1\\ 1 & 0\end{pmatrix},I_{i-1})A_1 & \text{if }1\le i\le n-1 \\
\diag_{(2n)\times(2n)}(I_{n-1},\begin{pmatrix} 0 & -1\\ 1 & 0\end{pmatrix},I_{n-1})A_1 & \text{if }i=n \\
\diag_{(2n)\times(2n)}(I_{2n-i-1},\begin{pmatrix} 0 & 1\\ 1 & 0\end{pmatrix},I_{2i-2n-2},\begin{pmatrix} 0 & 1\\ 1 & 0\end{pmatrix},I_{2n-i-1})A_1 & \text{if }n+1\le i\le 2n-1 \\
\end{cases}$$
share the same distribution. Set $M=(1+p^D\Mat_{(2n)\times(2n)}(\Z_p))\cap\GSp_{2n}(\Z_p)$, where $D>\mu_1^{(1)}-\mu_{2n}^{(1)}$ is sufficiently large. If $1\le i\le n-1$, let $S$ denote the set
\begin{multline*}
S=M\diag_{(2n)\times(2n)}(p^{\mu_2^{(1)}},\ldots,p^{\mu_i^{(1)}},\begin{pmatrix} p^{\mu_{i+1}^{(1)}} & 0\\ p^{\mu_{2n}^{(1)}} & p^{\mu_1^{(1)}}\end{pmatrix},p^{\mu_{i+2}^{(1)}},\ldots,p^{\mu_{2n-i-1}^{(1)}},\\
\begin{pmatrix} p^{\mu_{2n}^{(1)}} & 0\\ -p^{\mu_{2n-i}^{(1)}-\mu_1^{(1)}+\mu_{2n}^{(1)}} & p^{\mu_{2n-i}^{(1)}}\end{pmatrix},p^{\mu_{2n-i}^{(1)}},\ldots,p^{\mu_{2n-1}^{(1)}})M;\end{multline*}
If $i=n$, let $S$ denote the set
\begin{equation*}
S=M\diag_{(2n)\times(2n)}(p^{\mu_2^{(1)}},\ldots,p^{\mu_n^{(1)}},\begin{pmatrix} p^{\mu_{2n}^{(1)}} & 0\\ p^{\mu_{2n}^{(1)}} & p^{\mu_1^{(1)}}\end{pmatrix},p^{\mu_{n+1}^{(1)}},\ldots,p^{\mu_{2n-1}^{(1)}})M;
\end{equation*}
If $n+1\le i\le 2n-1$, let $S$ denote the set
\begin{multline*}
S=M\diag_{(2n)\times(2n)}(p^{\mu_2^{(1)}},\ldots,p^{\mu_{2n-i}^{(1)}},\begin{pmatrix} p^{\mu_{2n}^{(1)}} & 0\\ -p^{\mu_{2n-i+1}^{(1)}-\mu_1^{(1)}+\mu_{2n}^{(1)}} & p^{\mu_{2n-i+1}^{(1)}}\end{pmatrix},\\
p^{\mu_{2n-i+2}^{(1)}},\ldots,p^{\mu_{i-1}^{(1)}},\begin{pmatrix}p^{\mu_i^{(1)}} & 0\\ p^{\mu_{2n}^{(1)}} & p^{\mu_1^{(1)}}\end{pmatrix},p^{\mu_{i+1}^{(1)}},\ldots,p^{\mu_{2n-1}^{(1)}})M.\end{multline*}
By \Cref{prop:singular numbers of sum}, the matrices in the set $S$
have singular numbers $\mu^{(1)}=(\mu_1^{(1)},\ldots,\mu_{2n}^{(1)})$. The probability $\mathbf{P}(A_1\in S)>0$ is nonzero, and whenever this happens, \eqref{eq: i.i.d_inequality of weight for GSp} holds strict inequality. 
Therefore, taking expectations on both sides of \eqref{eq: i.i.d_inequality of weight for GSp} brings us the proof.
\end{proof}

Here are some possible ways that might lead to further research.

\begin{enumerate}
\item The universality results for $\GL_n$ (\Cref{cor: i.i.d.}) and $\GSp_{2n}$ (\Cref{thm: i.i.d for GSp}) demonstrate the importance of matrix corners in analyzing matrix products. Therefore, a natural direction for further study is to derive an explicit formula for the joint distribution of corners of a random matrix $A$, namely
$$\mathbf{P}(\SN(A^{(1)})=\mu^{(1)},\SN(A^{(2)})=\mu^{(2)},\ldots,\SN(A^{(i)})=\mu^{(i)},\ldots)$$
where $A^{(i)}$ is the submatrix of $A$, removing the first $i-1$ rows. For the $\GL_n$ case, such a distribution could be expressed through skew Hall-Littlewood polynomials (see \Cref{cor: joint distribution of corner}). Currently, such a precise form does not exist for the $\GSp_{2n}$ case, which we believe is worth studying and should be based on the form of Hall-Littlewood polynomials of type $C_n$ (see \cite{macdonald1990orthogonal} for constructions and details) because this is the Weyl group of the symplectic Lie algebra. 

\item The classification of root systems is well established in representation theory, see \cite[Chapter 2.4]{humphreys1992reflection}. The root system $A_n$ corresponds to the reductive group $\GL_n$, while the root system $C_n$ is related to $\GSp_{2n}$. However, for other types of root systems, it’s not always possible to assume the representatives $d(\lambda)$ to be diagonal. For example, the root systems $B_n$ and $D_n$ are related to $\SO_{2n+1}$ and $\SO_{2n}$ respectively. They are more complex because the orbits $\SO_n(\Z_p)\backslash\SO_n(\Q_p)/\SO_n(\Z_p)$ are more challenging to handle. Outside of types 
$A,B,C,D$, reductive groups may not even have an elegant form. Thus, it would be worthwhile to develop heuristics that can yield meaningful insights into these systems and prove the asymtotics if feasible.
\end{enumerate}

\appendix

\section{Hall-Littlewood polynomials and random matrix corners}\label{HL and RMT}

In this appendix, we introduce Hall-Littlewood polynomials and their applications in $p$-adic random matrices. Based on the form of such polynomials, we provide an alternate proof of \Cref{cor: i.i.d.} by direct computation. Unlike the non-negative partitions in Macdonald's book \cite{Macdonald}, we define Hall-Littlewood polynomials for finite-length integer signatures, which may include negative parts. This is because we are investigating $n\times n$ matrices and their $(n-i+1)\times n$ corners, where the length of the singular number tuple is naturally fixed. Furthermore, negative singular numbers must be considered because the entries are over $\Q_p$. For more background on Hall-Littlewood polynomials, please refer to \cite[Section 3]{Macdonald}.

\begin{defi}
Let $\lambda=(\l_1,\ldots,\l_n)\in\Sig_n$ be an integer signature of length $n$. Then, the \textbf{Hall-Littlewood polynomial} $P_\lambda(x_1,\ldots,x_n;t)$ is defined by

\begin{equation}\label{eq: Hall Littlewood polynomial}
P_\l(x_1,\ldots,x_n;t) = \frac{1}{v_\l(t)} \sum_{\sigma \in S_n} \sigma\left(x_1^{\lambda_1}\cdots x_n^{\lambda_n} \prod_{1 \leq i < j \leq n} \frac{x_i-tx_j}{x_i-x_j}\right)
\end{equation}
where $v_\l(t) = \prod_{i \in \Z} \frac{(1-t)\cdots(1-t^{m_i(\lambda)})}{(1-t)^{m_i(\l)}}$. 
\end{defi}

\begin{prop}
The Hall-Littlewood polynomials $P_\l(x_1,\cdots,x_n;t)$ satisfy the following properties:

\begin{enumerate}[left=0pt]
\item They have the form
$$P_\l(x_1,\ldots,x_n;t) = x_1^{\l_1}x_2^{\l_2}\cdots x_n^{\l_n} + \text{(lower-order monomials in the lexicographic order)};$$

\item When $\lambda\in\Sig_n^+$ ranges over all nonnegative integer signatures of length $n$, the Hall-Littlewood polynomials $P_\l(x_1,\ldots,x_n;t)$ form a $\Z[t]$ basis of $\L_n[t]$, where $\L_n[t]:=\Z[t][x_1,\ldots,x_n]^{S_n}$ is the ring of symmetric polynomials in $n$ variables $x_1,\ldots,x_n$ with coefficients in $\Z[t]$;

\item When $\lambda\in\Sig_n$ ranges over all integer signatures of length $n$, the Hall-Littlewood polynomials $P_\l(x_1,\ldots,x_n;t)$ form a $\Z[t]$ basis of $\Z[t][x_1^{\pm 1},\ldots,x_n^{\pm 1}]^{S_n}$, the ring of symmetric Laurent polynomials in $n$ variables $x_1,\ldots,x_n$ with coefficients in $\Z[t]$.
\end{enumerate}
\end{prop}

Because the $P_\l,\lambda\in\Sig_n$ form a basis for the vector space of symmetric Laurent polynomials in $n$ variables, there exist symmetric Laurent polynomials $P_{\l/\mu}(x_1,\ldots,x_{n-k};t) \in \Z[t][x_1^{\pm 1},\ldots,x_{n-k}^{\pm 1}]^{S_{n-k}}$ indexed by $\l \in \Sig_n, \mu \in \Sig_k$ which are defined by
\begin{equation*}\label{eq: skew Hall Littlewood polynomial}
    P_\l(x_1,\ldots,x_n;t) = \sum_{\mu \in \Sig_k} P_{\l/\mu}(x_{k+1},\ldots,x_n;t) P_\mu(x_1,\ldots,x_k;t).
\end{equation*}
where the polynomial $P_{\l/\mu}$ is named the \textbf{skew Hall-Littlewood polynomial}. The following lemma is deduced from (5.11'), Chapter 5 of Macdonald's book \cite{Macdonald}, which gives an explicit form of the skew Hall-Littlewood polynomial: 
\begin{lemma}
(Branching rule) For $\l=(\l_1,\ldots,\l_n)\in \Sig_n, \mu=(\mu_1,\ldots,\mu_{n-1}) \in \Sig_{n-1}$ with $\mu \prec_P \l$, let
\begin{equation}
    \psi_{\l/\mu}(t):= \prod_{\substack{i \in \Z \\ m_i(\mu) = m_i(\l)+1}} (1-t^{m_i(\mu)}).
\end{equation}
Then for $\l\in\Sig_n,\mu\in\Sig_{n-k}$, we have
\begin{equation}\label{eq: branching rule}
P_{\l/\mu}(x_1,\ldots,x_k;t)=\sum_{\mu = \l^{(1)} \prec_P \l^{(2)} \prec_P \cdots \prec_P \l^{(k)}= \l} \prod_{i=1}^{k-1} x_i^{|\l^{(i+1)}|-|\l^{(i)}|}\psi_{\l^{(i+1)}/\l^{(i)}}(t).
\end{equation} 
\end{lemma}

For the sake of convenience, we often write \eqref{eq: branching rule} in the abbreviated form
\begin{equation}\label{eq: abbreviate branching rule}
P_{\l/\mu}(x_1,\ldots,x_k;t) = \sum_{T \in \GT_P(\l/\mu)} \psi(T) \bx^{wt(T)}
\end{equation}
where $\GT_P(\l/\mu)$ is the set of sequences of interlacing signatures $\mu = \l^{(1)} \prec_P \l^{(2)} \prec_P \cdots \prec_P \l^{(k)}= \l$, $\bx$ is the collection of variables $x_1,\ldots,x_n$, and $wt(T) := (|\l^{(1)}|,|\l^{(2)}|-|\l^{(1)}|,\ldots,|\l^{(k)}|-|\l^{(k-1)}|) \in \Z^k$ is the \textbf{weight} of $T$.

We also give the value of the Hall-Littlewood polynomial under principal specialization, which serves as a corollary of \eqref{eq: Hall Littlewood polynomial}: for any $\lambda=(\lambda_1,\ldots,\lambda_n)\in\Sig_n$, we have
\begin{equation}\label{eq: principal specialization}
P_\l(x,xt,\ldots,xt^{n-1};t)= x^{|\l|} t^{n(\l)} \frac{(1-t)(1-t^2)\cdots(1-t^n)}{\prod_{i \in \Z}\prod_{j=1}^{m_i(\lambda)}(1-t^j)}.\end{equation}

The following theorem is essentially contained in \cite[Theorem 1.3]{Roger}, which is widely involved in our method:

\begin{thm}\label{thm: corner process and product process}
Set $t=1/p$.
\begin{enumerate}
\item Let $i,n$ be integers with $1 \le i  \le n$, $\mu^{(i)} \in \Sig_{n-i+1}$, and $A^{(i)}\in\Mat_{(n-i+1)\times n}(\Q_p)$ be random with $\SN(A^{(i)})=\mu^{(i)}$ fixed and distribution invariant under $\GL_{n-i+1}(\Z_p) \times \GL_n(\Z_p)$ acting on the left and right. Let $A^{(i+1)} \in\Mat_{(n-i) \times n}(\Q_p)$ be the last $n-i$ rows of $A$. Then $\SN(A^{(i+1)})$ is a random element of $\Sig_{n-i}$ with distribution
\begin{equation}\label{eq:p-adic_corners}
\Pr(\SN(A^{(i+1)}) = \mu^{(i+1)}) = P_{\mu^{(i)}/\mu^{(i+1)}}(1;t) \frac{P_{\mu^{(i+1)}}(t,\ldots,t^{n-i};t)}{P_{\mu^{(i)}}(1,\ldots,t^{n-i};t)}.
\end{equation} 

\item Let $A\in\Mat_{n\times m}(\Z_p)$ be the same as \Cref{ex: nxm corner of invertible}. Then the distribution of $\SN(A)$ is given by the Hall-Littlewood measure
$$\mathbf{P}(\SN(A)=\lambda)=\frac{P_\lambda(1,t,\ldots,t^{n-1};t)Q_\lambda(t^{m-n-1},\ldots,t^{N-n};t)}{\Pi_t(1,t,\ldots,t^{n-1};t^{m-n-1},\ldots,t^{N-n})}$$
Here $Q_\lambda$ is the Hall-Littlewood $Q$ polynomial, and $\Pi_t$ is the Cauchy Kernel, see \cite[Section 3]{Macdonald}.
\end{enumerate}
\end{thm}

As a corollary of the above, we iterate \eqref{eq:p-adic_corners} and get the following:

\begin{cor}\label{cor: joint distribution of corner}

Let $A\in\Mat_{n\times n}(\Q_p)$ be random with $\SN(A)=\mu^{(1)}$ fixed and distribution invariant under the action of $\GL_n(\Z_p)\times\GL_n(\Z_p)$ acting on the right and left. For every $1\le i\le n$, let $A^{(i)}$ be the submatrices of $A$, consisting of the last $n-i+1$ rows and all $n$ columns. Then the joint distribution of $(\SN(A^{(2)}),\ldots,\SN(A^{(n)}))$ has the form 
$$\mathbf{P}(\SN(A^{(2)})=\mu^{(2)},\ldots,\SN(A^{(n)})=\mu^{(n)})=\frac{P_{\mu^{(1)}/\mu^{(2)}}(1;t)\cdots P_{\mu^{(n-1)}/\mu^{(n)}}(t^{n-2};t)P_{\mu^{(n)}}(t^{n-1};t)}{P_{\mu^{(1)}}(1,t,\ldots,t^{n-1};t)}.$$

In particular, for all $2\le k\le n$, the distribution of $\SN(A^{(k)})$ has the form
$$\mathbf{P}(\SN(A^{(k)})=\mu^{(k)})=\frac{P_{\mu^{(1)}/\mu^{(k)}}(1,t,\ldots,t^{k-2};t)P_{\mu^{(k)}}(t^{k-1},t^k,\ldots,t^{n-1};t)}{P_{\mu^{(1)}}(1,t,\ldots,t^{n-1};t)}.$$
\end{cor}

\begin{rmk}\label{rmk: Roger setting split}
Now we return to the discussion in \eqref{item: Roger setting} of \Cref{example: split}. According to \Cref{thm: corner process and product process}, the joint distribution of $(\SN(A_k^{(1)}),\ldots,\SN(A_k^{(n)}))$ takes the form 
$$\mathbf{P}(\SN(A_k^{(1)})=\mu^{(1)},\ldots,\SN(A_k^{(n)})=\mu^{(n)})=\frac{Q_{\mu^{(1)}}(t,\ldots,t^{N_k-n};t)P_{\mu^{(1)}/\mu^{(2)}}(1;t)\cdots P_{\mu^{(n)}}(t^{n-1};t)}{\prod_t(1,t,\ldots,t^{n-1};t,\ldots,t^{N_k-n})}.$$

Recall the property of the Hall-Littlewood process (see Proposition 2.2.6 of \cite{borodin2014macdonald}), the random variables $|\SN(A_k^{(1)})|-|\SN(A_k^{(2)})|,\ldots,|\SN(A_k^{(n-1)})|-|\SN(A_k^{(n)})|,|\SN(A_k^{(n)})|$ are independent, and we have

$$\E(x^{|\SN(A_k^{(j)})|-|\SN(A_k^{(j+1)})|})=\frac{\prod_t(x;t^j,\ldots,t^{N_k-n+j-1})}{\prod_t(1;t^j,\ldots,t^{N_k-n+j-1})},\quad\forall 1\le j\le n-1$$

$$\E(x^{|\SN(A_k^{(n)})|})=\frac{\prod_t(x;t^n,\ldots,t^{N_k-1})}{\prod_t(1;t^n,\ldots,t^{N_k-1})}$$
This coincides with definition 12 of \cite{Roger}(In other words, the sequence $v$ we defined has the same distribution as in the previous paper.). Hence, the assertion that the sequence $v$ splits almost surely has already been verified by \cite[Lemma 5.5]{Roger}.
\end{rmk}

In the end, we give another proof of \Cref{cor: i.i.d.}, based on the joint distribution given by \Cref{cor: joint distribution of corner}.
\begin{proof}
Suppose there exists $\mu^{(1)}\in\Sig_n$ such that $\mathbf{P}_1(\mu^{(1)}):=\mathbf{P}(\SN(A_1)=\mu^{(1)})>0$ and $\mu_1^{(1)}>\mu_n^{(1)}$. We only need to prove the inequality \eqref{eq: inequality of corners}. When we compute the expectations $\E|\SN(A_1^{(1)})|,\ldots,\E|\SN(A_1^{(n)})|$, we first compute the case that $\SN(A_1)=\mu^{(1)}$ is fixed, then add these results according to the distribution $\mathbf{P}_1$. Hence, there is no loss we study the case $\mathbf{P}_1(\mu^{(1)})=1$ is fixed. By \Cref{cor: joint distribution of corner}, we have

$$\mathbf{P}(\SN(A_1^{(j)})=\mu^{(j)})=\frac{P_{\mu^{(1)}/\mu^{(j)}}(1,t,\ldots,t^{j-2};t)P_{\mu^{(j)}}(t^{j-1},t^j,\ldots,t^{n-1};t)}{P_{\mu^{(1)}}(1,t,\ldots,t^{n-1};t)}.$$
Therefore, for all $x\in\R$, we have
\begin{align}
\E(x^{|\SN(A_1^{(j)})|})&=\sum_{\mu^{(j)}}x^{|\mu^{(j)}|}\frac{P_{\mu^{(1)}/\mu^{(j)}}(1,t,\ldots,t^{j-2};t)P_{\mu^{(j)}}(t^{j-1},t^j,\ldots,t^{n-1};t)}{P_{\mu^{(1)}}(1,t,\ldots,t^{n-1};t)}\nonumber\\
&=\sum_{\mu^{(j)}}\frac{P_{\mu^{(1)}/\mu^{(j)}}(1,t,\ldots,t^{j-2};t)P_{\mu^{(j)}}(xt^{j-1},xt^j,\ldots,xt^{n-1};t)}{P_{\mu^{(1)}}(1,t,\ldots,t^{n-1};t)}\nonumber\\
&=\frac{P_{\mu^{(1)}}(1,\ldots,t^{j-2},xt^{j-1},\ldots,xt^{n-1};t)}{P_{\mu^{(1)}}(1,t,\ldots,t^{n-1};t)}.
\end{align}
Taking the derivative on both sides and set $x=1$, we obtain
$$\E|\SN(A_1^{(j)})|=\frac{1}{P_{\mu^{(1)}}(1,t,\ldots,t^{n-1};t)}\frac{d}{dx}P_{\mu^{(1)}}(1,\ldots,t^{j-2},xt^{j-1},\ldots,xt^{n-1};t)\Big|_{x=1},\quad\forall 1\le j\le n.$$
Therefore, for all $1\le j\le n-1$, the difference $\E|\SN(A_1^{(j)})|-\E|\SN(A_1^{(j+1)})|$ has the form
$$\E|\SN(A_1^{(j)})|-\E|\SN(A_1^{(j+1)})|=\frac{1}{P_{\mu^{(1)}}(1,\ldots,t^{n-1};t)}\frac{d}{dx}P_{\mu^{(1)}}(1,\ldots,t^{j-2},xt^{j-1},t^j,\ldots,t^{n-1};t)\Big|_{x=1}$$

The following lemma is needed:

\begin{lemma}
For any $\mu^{(1)}=(\mu_1^{(1)},\ldots,\mu_n^{(1)})\in\Sig_n$, let $m_{\mu^{(1)}}(x_1,\ldots x_n)=\sum_\alpha x^\alpha$, summed over all distinct permutations $\alpha=(\alpha_1,\ldots,\alpha_n)$ of $\mu^{(1)}$ (See Subsection 1.2 of \cite{Macdonald}). Let $b_1>b_2>\ldots>b_n>0$. Then we have
$$\frac{d}{dx}m_{\mu^{(1)}}(xb_1,b_2,\ldots,b_n)\Big|_{x=1}\ge\ldots\ge\frac{d}{dx}m_{\mu^{(1)}}(b_1,b_2,\ldots,xb_n)\Big|_{x=1}.$$
Every equality holds if and only if $\mu_1^{(1)}=\ldots=\mu_n^{(1)}$.
\end{lemma}

\begin{proof}
If $\mu_1^{(1)}=\ldots=\mu_n^{(1)}$, then $m_{\mu^{(1)}}(xb_1,b_2,\ldots,b_n)=\ldots=m_{\mu^{(1)}}(b_1,b_2,\ldots,xb_n)=(xb_1\ldots b_n)^{\mu_1^{(1)}}$, and the equalities hold. Now suppose $\mu_1^{(1)}>\mu_n^{(1)}$. We only prove the inequality 
$$\frac{d}{dx}m_{\mu^{(1)}}(xb_1,b_2,\ldots,b_n)\Big|_{x=1}>\frac{d}{dx}m_{\mu^{(1)}}(b_1,xb_2,\ldots,b_n)\Big|_{x=1},$$ 
and the other cases are the same. In fact, we have
$$\frac{d}{dx}m_{\mu^{(1)}}(xb_1,b_2,\ldots,b_n)\Big|_{x=1}=\sum_{\alpha}\alpha_1b_1^{\alpha_1}\ldots b_n^{\alpha_n},$$
$$\frac{d}{dx}m_{\mu^{(1)}}(b_1,xb_2,\ldots,b_n)\Big|_{x=1}=\sum_{\alpha}\alpha_2b_1^{\alpha_1}\ldots b_n^{\alpha_n}.$$
Therefore, their difference has the form
$$
\frac{d}{dx}(m_{\mu^{(1)}}(xb_1,b_2,\ldots,b_n)-m_{\mu^{(1)}}(b_1,xb_2,\ldots,b_n))\Big|_{x=1}=\sum_{\alpha}(\alpha_1-\alpha_2)b_1^{\alpha_1}\cdots b_n^{\alpha_n}$$
$$=\frac{1}{2}\sum_{\alpha}(\alpha_1-\alpha_2)(b_1^{\alpha_1}b_2^{\alpha_2}-b_1^{\alpha_2}b_2^{\alpha_1})b_3^{\alpha_3}\cdots b_n^{\alpha_n}
$$
which must be positive since $\mu_1^{(1)}>\mu_n
^{(1)}>0$. This ends the proof of the lemma.
\end{proof}
Turning back to the proof, according to the branching rule given by \eqref{eq: branching rule}, the Hall-Littlewood polynomial can be expressed by the sum of monomials $P_{\mu^{(1)}}=\sum_\nu c_{\mu^{(1)}\nu}m_\nu$, where $c_{\mu^{(1)}\nu}=c_{\mu^{(1)}\nu}(t)$ is non-negative, and equals one when $\nu=\mu^{(1)}$. Apply the above lemma, and set $b_1=1,b_2=t,\ldots,b_n=t^{n-1}$, we end the proof by proving the inequality \eqref{eq: inequality of corners} holds.
\end{proof}


\end{document}